\font\steptwo=cmb10 scaled\magstep2

\font\tt=cmtt10 \relax
\magnification=\magstep1
\settabs 18 \columns

\hsize=16truecm

\def\b{\bigskip}
\def\bb{\bigskip\bigskip}

\def\no{\noindent}

\def\ce{\centerline}
\def\ve{\vfill\eject}

\def\cb{\hbox{${\cal B}$}}

\def\harr#1#2{\smash{\mathop{\hbox to .25 in{\rightarrowfill}}
 \limits^{\scriptstyle#1}_{\scriptstyle#2}}}

\def\today{\ifcase\month\or January\or February\or March\or April\or  May\or
June\or July\or August\or September\or October\or November\or  December\fi
\space\number\day, \number\year }

\def\ui{{\underline i}}

\def\p{\partial}

\def\today{\ifcase\month\or January\or February\or March\or April\or  May\or
June\or July\or August\or September\or October\or November\or  December\fi
\space\number\day, \number\year }

 \b

\ce{{\steptwo {ON THE CLASSIFICATION OF q-ALGEBRAS}}}
\b
\ce{Christian Fr\o nsdal}
\b
\ce{\it Physics Department, University of California, Los Angeles CA 90095-1547
USA}

\bb\bb

\no{\it ABSTRACT.} The problem is the classification of the ideals of ``free
differential algebras", or the associated quotient algebras, the  q-algebras;
being finitely generated,  unital ${\bf C}$-algebras with homogeneous relations
and a q-differential structure. This family of algebras includes the quantum
groups, or at least those that are based on simple (super)  Lie or Kac-Moody
algebras. Their classification would encompass the so far incompleted 
classification of quantized (super) Kac-Moody algebras and of the (super)  
Kac-Moody algebras themselves. These can be defined as singular limits of 
$q$-algebras, and it is evident that to deal  with the $q$-algebras in their 
full generality is more rational  than the examination of each singular limit
separately. This is not just because quantization unifies algebras and
superalgebras, but also because the points  ``$q = 1$" and
$``q = -1$" are the most singular points in parameter space. In this paper one
of two major hurdles in this classification program has been overcome.  Fix a
set  of integers $n_1,...,n_k$, and consider the space $\cb_Q$  of 
homogeneous polynomials of degree $n_1$ in the generator $e_1$, and so on.  
Assume that there are no constants among the polynomials of lower degree, 
in any one of the generators; in this case all constants in the space 
$\cb_Q$ have been classified. The task that remains, the more formidable
one, is to remove the stipulation that there are no constants of lower degree.

\b

\no{\bf 1. Introduction.}

The classification of simple Lie algebras, achieved by Killing and Cartan near
the end of the 19'th century, was a milestone in modern mathematics. Other
classification problems came later: simple Lie bialgebras, super Lie algebras
[BD], Kac-Moody algebras [Kc], [M], generalized Kac-Moody algebras, quantum
groups. Little is known about generalized Kac-Moody algebras, and the
classification of Lie superalgebras is only now nearing completion.  In fact, a
complete understanding of the  relations that govern the latter turned out to
be a quite formidable task  [MZ],
[GL], [B], [FSS] . It was achieved with the methods of  combinatorial
algebras, Groebner bases and Lyndon words, tools the effectiveness of which
ultimately relies on  the Diamond Lemma [B]. We are reminded of
Bergman's remark, that  some problems become more tractable  when generalized
to an  ultimate degree. It may be that some of the struggle associated with
super Lie algebras can be reduced by generalization. From the point of view of
quantum groups, Lie algebras, superalgebras and Kac-Moody algebras are very
singular limits. This paper suggests that the difficulties can be mitigated by
working with quantum groups and passing to the interesting limits afterwards.

We  deal with a class of algebras that forms a natural generalization
of the quantum groups first defined systematically by Drinfel'd [D]. Simplicity
is gained by a method first suggested by Lusztig [L]. In this formulation the
Cartan subalgebra disappears and there remains only the subalgebra generated by
the positive simple root vectors. This  structure also appears in the
work of  Kashiwara [Ks] and others in connection with crystal
bases.   The trick of eliminating the Cartan subalgebra from the
scene was rediscovered by the author in his work on universal R-matrices [F].
It led to the study of ``free differential algebras", for which a better term
might be $q$-algebras.

One begins with the free, unital ${\bf C}$-algebra \cb ~on generators
$\{e_i\}_{i \in {\cal N}  = \{1,...,N\}}$ that should be regarded as
positive Serre generators. The negative Serre generators are represented
(replaced) by q-differential operators. The {\it parameters} are the values of
a function $q: {\cal N} \times {\cal N} \rightarrow {\bf C},~ (i,j) \mapsto
q_{ij}
\neq 0$. They appear in the action of a set $\{\partial_i\}_{i \in {\cal N}}$
of q-differential operators on \cb, this action being defined
by $\p_i x = 0 $ for $x \in {\bf C} \subset \cb$ and  iteratively for all
$x \in \cb$
$$
\p_i(e_jx) = \delta_i^j x + q_{ij}e_j\partial_i x\,.\eqno(1.1)
$$ Henceforth \cb~ will stand for the algebra \cb~ endowed with this
$q$-differential structure. It is naturally graded, $\cb = \bigoplus_{n =
0,1,...} \cb_n$, with $\cb_0 = {\bf C}$. A {\it constant} is an element
$C\in \bigoplus_{n = 1,2,...}\cb_n$ with the property $\p_i C = 0,~ i =
1,...,N$. For parameters $\{q_{ij}\}$ in general position there are no
constants, constants exist only for values of the parameters that by that
virtue will be called {\it singular}. The constants generate an ideal
$I_q\subset \cb$. The algebras  of interest (quantum groups, quantized
Kac-Moody algebras among them, see below) are the quotient algebras $\cb' =
\cb/{\cal I}_q$, with the inherited q-differential structure. Here are some
examples. In each case the parameters are supposed to be generic except for the
explicitly imposed constraints.
\b (1) If for all $i \in {\cal N}$ we have $q_{ii} = -1$, then $I_q$ is
generated by $\{e_i^2\}_{i \in {\cal N}}$. (2) If for all $(i,j),~i \neq j$ we
have $q_{ij}q_{ji} = 1$, then $\cb'$ is Manin's quantum plane. (3) If for all
$(i,j),~i\neq j$ there is a non-negative integer $k_{ij}$ such that
$$ (q_{ii})^{k_{ij}}q_{ij}q_{ji} = 1,\eqno(1.2)
$$ then $\cb'$ is a quantized, generalized Kac-Moody algebra with Cartan matrix
$A_{ij} = -k_{ij}$ for $i \neq j$. This last category includes quantized
versions of the basic, generalized super-Kac-Moody algebras; that is, those
that are based on a Cartan matrix. In fact, the distinction between Lie
algebras and super Lie algebras disappears upon quantization (in the sense of
Drinfel'd) and this may be taken as an indication that certain aspects of the
theory actually become simpler after quantization. (4) A final example
[F],[FG],
to show that this category of algebras encompasses new and interesting
possibilities,  is the case that the parameters satisfy the single constraint
$\prod q_{ij} = 1$, where the product runs over all imbeddings of $\{i,j\}$
into $\{1,...,n\}$. The ideal is generated by a polynomial in $N$ variables.
Examples of this kind play a role in quantized super Lie
algebras\footnote*{Thanks to Cyrill Oseledetz for pointing this out.} [FSS].

The problem of classification that we are working on here is a very natural
one:

\hskip1cm {\it Find all the singular points in parameter

\hskip1cm  space and describe the associated ideals in \cb.}

 Of course we shall not accomplish this. But we shall complete
 the following program. First we remark that the space of constants is
generated by homogeneous constants. Let $Q$ be a set of natural numbers, with
repetition, with cardinality $n$, and  let $\cb_Q$ be the subspace of \cb
~generated by the set of monomials
$\{e_{i_1}...e_{i_n}\}$ where the indices run over all permutations of
$Q$. Program:

\hskip1cm{\it Under the stipulation that there are no constants in $\cb_{Q'}$,

\hskip1cm
$Q'$ any proper subset of $Q$, find all the singular points,

\hskip1cm  and all the associated constants, in $\cb_Q$}.

This problem had already been solved for the case $Q = \{1,2,...,n\}$. The
complications introduced by repetitions in
$Q$ will be recognized as similar to those that arise when one passes from Lie
algebras to super Lie algebras and color super Lie algebras [MZ], and the
methods of analysis will overlap. We find that the introduction of a Groebner
basis in terms of generalized Lyndon words is useful, but  another basis has
been more effective, so far. The answer is given in terms of certain
idempotents
in the group algebra of $S_n$, and the lengths of certain orbits of cyclic
subgroups acting in $S_nQ$.
\b

\no{\it Relation to other work.}

The problematics of this paper is of course related to the vast litterature on
quantum groups, especially to references already quoted. However, it is
essentially different in that it explores a larger category. (We use the word
in a non-technical sense.) It is sufficiently larger than the usual one of
quantum groups to be interesting. There is the hope that some
results may be easier to come by in this general setting. To my knowledge, this
point of view is found in the following papers. First, there is the work of
Varchenko [V] on the arrangement of hyperplanes. His emphasis is on quantum
groups, but the general situation is clearly envisaged. The paper [R]
explicitly calls for an exploration of generalizations of quantum groups, in
precisely the same direction. But the work of Kharchenko [Kh] appears to be
closest in scope to our work. This is perhaps not obvious at first sight, so it
is necessary to sketch the Hopf algebra aspects of the algebras that are under
investigation.

Consider the automorphisms $\{{\cal K}_i,K^i\}_{i = 1,...,N}$ defined by
$ {\cal K}_i:  e_j \mapsto (q_{ji})^{-1}e_j$, ${\cal K}^i:  e_j \mapsto
q_{ij}e_j$.
  Enlarge the algebra \cb~ by adjoining invertible elements $K_i, K^i$ that
implement these automorphisms, and consider the algebra generated by ${e_i,
K_i, K^i}$ with relations
$$
K_ie_j(K_i)^{-1} = (q_{ij})^{-1}e_j.
$$
This becomes a bialgebra  with coproduct and antipode defined by
$$
\Delta e_i = e_i \otimes 1 + K_i \otimes e_i,~~\Delta K_i = K_i \otimes
K_i,~~ i = 1,...,N.
$$
$$
S(K_i) = (K_i)^{-1},~ S(e_i) = - K_ie_i.
$$
Kharchenko studies certain
``quantum relations" in this algebra; they are precisely the constants. In the
paper [Kh], is a result on the classification of these relations that coincides
with a result in [FG].  The relationship to Universal R-matrices and Kac-Moody
algebras was summarized in [FG], but only the paper [F] goes into detail.

\b
\no{\it Summary and results. }

\no There is an old, partial result. Its scope limited in two ways: {

 \narrower   { (1) It is assumed that the parameters are such that there are no
constants in $\cb_{Q_j}$, $Q_j$ the  subset of
$Q$  obtained from $Q$ by reducing by 1 the incidence of $j$ in $Q$. In other
words, one looks for new phenomena that first appear at polynomial order $n$ in
the ``most generic" case. }

(2)  The degrees of homogeneity are restricted to 0 and 1. That is, one
considers, for some natural number $n$, the subspace $\cb_{\underline n}\,,~
\underline n = 12...n\,$,  of homogeneous polynomials of degree one in each
generator
$e_1,...,e_n$.\smallskip }
\no Within the limitation of both stipulations, it was possible to get complete
results [FG]. (See also [V] and [Kh].) The space of constants in
$\cb_{\underline n}$ is empty unless
$$
\prod q_{ij} = 1,\eqno(1.3)
$$ and then  the dimension of the space of constants is $(n-2)!$.  An algorithm
was given for the construction of all the constants in this case.

In this paper only the second of the two stipulations will be removed. We
consider polynomials spanned by $e_{i_1}...e_{i_n}$, with fixed degrees of
homogenity, $\underline i = \{i_1...i_n\}$ running over the permutations $\hat
Q$ of an arbitrary set
$Q = 1^{n_1}2^{n_2}...,~ \sum n_i = n$.   To state the result we must consider
the natural action of $S_n$ in $\hat Q$ and in $\cb_Q$, as well as the action
of $S_{n-1}$ in $\hat Q_j$ and in $\cb_{Q_j}$. In the case $Q = \underline
n\,$, the product (1.3) runs over all pairs $(i,j) \subset \underline n\,$, in
general (1.3) is interpreted as a cocycle condition  and the product runs over
orbits of the cyclic subgroups
$S_n^c$ and $S^c_{n-1}$. Complications arise whenever some of these orbits are
short.

{\narrower \noindent{}{For a choice of the parameters $q_{ij}$, let
$\chi,~\chi_j$ be the number of orbits of $S_n^c$ in $\hat Q$, resp.
$S_{n-1}^c$ in $\hat Q_j$, for which the cocycle condition (1.3) is verified,
then {\it the dimension of the space of constants is} $\sum \chi_j - \chi$.
Define certain iterated $q$-commutators $X^{\ui},~\ui = i_1...i_n$~, see
Section 2.1.3, with the property that
$$
\p_j X^{\ui} = 0, ~j \neq i_1,
$$ and let $\hat X^j = \{X^\ui\,;~ i_1 = j\}$. {\it The subspace spanned by any
one of these sets of polynomials contains all the constants.} This gives a
practical method for the computation of the constants.\smallskip}}

\no The result is both complete and of a comforting simplicity, but the
demonstration is not satisfying, for we were unable to base the proof on a
coherent strategy. Neither the old result, for the case $Q = \underline n$, nor
the method used to obtain it,  was brought to bear on the more general problem,
because that method was lacking in elegance and probably not suited for more
complicated situations. We first examine a subspace spanned by simple
commutators of the form
$[e_{i_1}...e_{i_{n-1}},e_{i_n}]_{a(\underline i)}$, with the commutation
factor $a$ satisfying a certain cocycle condition. This space is certainly
fundamental;  it is the same as that spanned by the iterated commutators
$X^\ui$ and  it was proved that it contains all the constants. The relation to
Kac-Moody algebras and super Kac-Moody algebras is evident and it is felt that
these q-commutators are easier to deal with than the (anti-) commutators that
arise in the limit  when the commutation factors tend to
$\pm 1$. This belief is based on the observation that those limiting cases are
the most singular points in parameter space. To get this far we also had to
find a trick that relates the constants in
$\cb_Q$ to the constants in $\cb_{\underline n}$. A third method had to be used
to finish the job.   A basis built on Lyndon words, and an identification of
Lyndon words with iterated commutators, of the type used in connection with
superalgebras [MZ], gave us an existence lemma that was required to reach our
goal. But that is the only effective use that we were able to make of this
basis.

The first stipulation was essential. It is unclear which, if  any,
of the several methods used here will turn out to be of use in the future, when
we shall try to eliminate it.

\b
\no {\bf 2. Constants and commutators.}

 For any ordered set $n_1,...,n_k$ of positive integers  set
$$ n: = \sum_{i = 1}^kn_i,
$$  and let $Q$ be the ordered set,
$$ Q = \{1,...,1,2,...,2,...,k,...,k\}= \{(1)_{n_1},(2)_{n_2},...,(k)_{n_k}\}.
$$ Let $\cb_Q$ be the linear space spanned by  the polynomials
$e_{i_1},...,e_{i_n}$, where $i_1,...,i_n$ runs over the set
$\hat Q := S_nQ$ of distinct permutations of $Q$.  This paper studies the
constants in $\cb_Q$ in the case that there are no constants in
$\cb_{Q'}$, for any proper subset $Q'\subset Q$.
\b
\no{\it 2.1. Simple and iterated commutators..}

 \no{\it 2.1.1. Simple commutators and cyclic permutations.}  We write $\ui$
for $i_1...i_n$, and sometimes for other index sets, when the context makes it
clear which set is meant.  Fix a function
$a:\hat Q\rightarrow {\bf C}-\{0\},~\ui
\mapsto  a(\ui) \neq 0$ and let $\hat A_Q$ be the set of simple commutators
$$
\hat A_Q := \{[e_{i_1}...e_{i_{n-1}},e_{i_n}]_{a(\ui)};~ \ui \in \hat
Q\}.\eqno(2.1)
$$ Let $S^c_n$ denote the cyclic subgroup of $S_n$, and $[s Q ]^c \subset \hat
Q$ the orbit of
$S_n^c$ through $s Q \in \hat Q$. The function $a$ may be called a {\it
commutation factor} ; compare [MZ].
\b
\no{\it 2.1.2. Proposition.} For $\{a(\ui)\}$ in general position, the set
$\hat A_Q$ spans
${\cal B}_Q$. The full set of linear relations in ${\bf C}\hat A_Q$ is
generated by relations of the following type. For each $s Q  \in \hat Q$, such
that
$$
\prod_{  \ui\,\in\, [  s Q ]^c} a( \ui) = 1\,;\eqno(2.2)
$$
 that is, such that $a(\ui),[s Q]^c$ defines an $S_n^c$ cocycle, there is a
unique relation
$$
\sum_{ \ui\,\in\, [s Q]^c}C(\ui)[e_{i_1}...e_{i_{n-1}}, e_{i_n}]_{a(\ui)} =0
 \,,\eqno(2.3)
$$ with complex coefficients $C(\ui),\, \ui \in [sQ]^c$ all different from
zero.
\b More precisely, this relation has the form
$$
\sum_{\alpha = 0,1,...}\biggl(\prod_{\beta <\alpha  }a(\tau^\beta
\ui)\biggr)~\tau^\alpha [e_{i_1}...e_{i_{n-1}},e_{i_n}]_{a(\ui)} = 0.\eqno(2.4)
$$ where $\ui = s Q$ ($s \in S_n$) is fixed,   $\tau$ is a generator of $S_n^c$
and the sum extends to the orbit of $S_n^c$ through $s Q$. The operator
$$ P(\ui) :=
\sum_{\alpha = 0,1,...}\biggl(\prod_{\beta = 0,...,\alpha -1}a(\tau^\beta
\ui)\biggr)~\tau^\alpha
\eqno(2.5)
$$   is thus, up to a numerical factor, an idempotent in the group algebra of
$S_n$. When
$\hat A_Q$ fails to span $\cb_Q$, a supplemental basis is given by any  set of
monomials that contains exactly one element from each orbit of $S_n^c$ for
which the cocycle condition holds. More generally,  the {\it
orbital polynomial}
$$ u = \sum_\alpha C_\alpha \tau^\alpha e_{i_1}...e_{i_n}\eqno(2.6)
$$ is in ${\bf C}\hat A_Q$ if and only if it satisfies the {\it dual cocycle
condition}
$$
\sum_\alpha\biggl(\prod_{\beta \geq \alpha} a(\tau^\beta\ui)\biggr)C_\alpha =
0,\eqno(2.7)
$$ and a supplemental basis is made up of any collection of orbital
polynomials, one for each orbit, provided that they all violate the dual
cocycle condition.
\b
\no{\it 2.1.3. Iterated commutators.} For any permutation $sQ = \{i_1...i_n\}$
of $Q$,
$s\in S_n$, and commutation factors
$a(i_1...i_p) \neq 0,~ p = 1,...,n$,  let
$$ X^{i_1} = e_{i_1},~~X^{i_1...i_p} =
[X^{i_1...i_{p-1}},e_p]_{a(i_1...i_p)},~~ p = 1,...,n.
$$
\b
\no{\it 2.1.4. Proposition.} For commutation factors $a(i_1...i_p)$ in general
position the set
$$
\hat X :=\bigl\{X^{i_1...i_n}\,;~ i_1...i_n \in  \hat Q\}\eqno(2.8)
$$ spans $\cb_Q$.
\b
\no{\it Proof.} By induction; in fact, in this case ${\bf C}\hat X = {\bf
C}\hat A_Q$.
\b
\ve
\no{\it 2.2. Constants.}

\no{\it 2.2.1.} From now on the commutation factors will be chosen as follows.
Fix a function $q:
\hat Q
\otimes
\hat Q \rightarrow {\bf C}-\{0\},~ (i,j) \mapsto q_{ij} \neq 0$ and set
$$ a(i_1...i_p) = q_{i_pi_1}...q_{i_pi_{p-1}}.\eqno(2.9)
$$ The complex numbers $q := \{q_{ij}\}$ are {\it the
parameters}.\footnote*{For another interpretation see 3.2.6.} The action of the
$q$-differential operators $\p_i$ in \cb ~was defined in Section 1. Let $u =
e_{i_1}...e_{i_{p-1}}$, then for all $x\in\cb$,
$$
\p_{i_p}(ux) = (\p_{i_p}u)x + a(u,i_p)u\p_{i_p}x,~~a(u,i_p) := a(i_1...i_p),
$$  and if $p = n$, then $[u,e_{i_n}]_{a(u,i_n)} \in \hat A_Q$. The mapping
defined by $u \mapsto
\{a(u,i)\}_{i=1,...,N}$ is a grading of \cb; graded elements are said to be
{\it homogeneous}.

The action of $e_i$ on $\cb$ and on the quotient $\cb' = \cb/{\cal I}_q$
defined by $u \mapsto [u, e_i]_{a(u,i)}$  generalizes the adjoint
action of quantum groups. (A proof of this statement will be published.)

The significance of the choice (2.9)  in the present context  is that
 $$
\p_j[u,e_i]_{a(u,i)} = [\p_ju,e_i]_{a(u,i)q_{ji}},\eqno(2.10)
$$ for all pairs $(i,j)$ including the case that $i = j$, and that,
consequently,
$$
\p_j X^{i_1...i_p} = 0,~ j \neq i_1.\eqno(2.11)
$$
The cocycle condition (2.2) reduces, in the case of long orbits (orbits of
 length $n$) to
$$
\prod q_{ij} = 1,\eqno(2.12)
$$ where the product runs over all $n(n-1)$ pairs $(i,j) \subset \hat Q$.
Example: $Q = 1123$, Eq.(2.12) reads
$(q_{11}  q_{12} q_{21} q_{13} q_{31})^2 q_{23}q_{32} = 1$.
\b
\no{\it 2.2.2.} Let $Q_i, i = 1,...,k$ be the sets obtained from $Q$ by
removing one element, thus reducing $n_i$ by one, $i = 1,...,k$, respectively,
$\hat Q_j$ the set of permutations of $Q_j$ and $\cb_{Q_j}, j = 1,...,k$ the
corresponding spaces of polynomials.
\b
\no{\it 2.2.3. Stipulation  }.  It is supposed from now on, throughout Sections
2 and 3, that the parameters $q_{ij}$  are such that,     for each $j
= 1,...,k$ separately,  (a) $\hat A_{Q_j}$ spans $\cb_{Q_j}$ and (b)  there are
no constants in $\cb_{Q_j}\,$; for each $j = 1,...,k$.
\b The case $Q = \{1^n\}$ shows that (a) and (b) are independent. But it will
turn out that, if $Q \neq 1^n$, then (b) implies (a). See corollary 2.2.6.
 \b
 \ve
\no{\bf 2.2.4. Theorem.} If $C\in \cb_Q$ is a constant:
$\
\p_jC = 0$ for $j = 1,...,k$, then $C \in {\bf C}\hat A_Q$.
\b

\no {\it Proof.}  Let   $\{{\cal O}_\alpha\}$ be the collection of orbits of
$S_n^c$ in $\hat Q$ or in $\cb_Q$  on which the cocycle condition (2.2) holds.
In view of Eq.(2.9), this set is non-empty iff and only if
$$
\prod q_{i_1i_2} = 1.\eqno(2.13)
$$
 We refer to them as the {\it singular} orbits. Select a monomial
$u(\alpha) \neq 0$
  from each singular orbit ${\cal O}_\alpha$ of $S_n^c$ in
$\cb_Q$. A monomial cannot satisfy the cocycle condition, therefore
$$
\cb_Q = \bigoplus_\alpha{\bf C} u(\alpha) ~\oplus~ {\bf C}\hat A_Q.
\eqno(2.14)
$$
 A constant in $\cb_Q$ is a polynomial
$$ C = \sum_\alpha C(\alpha)\, u(\alpha) + \sum_{\ui\in\hat
Q}C'(\ui)[e_{i_1}...e_{i_{n-1}},e_{i_n}]_{a (\ui)},\eqno(2.15)
$$ such that, for $j = 1,...,k$,
$$
\p_jC= \sum_\alpha C(\alpha)\, \p_ju(\alpha) +
\sum_{\ui\in\hat
Q}C'(\ui)[\p_j(e_{i_1}...e_{i_{n-1}}),e_{i_n}]_{a(\ui)q_{ji_n}}  =
0.\eqno(2.16)
$$ Eq.(2.10) was used. The existence of a nontrivial solution of (2.16) implies
the existence of a nontrivial solution of the following equation,
$$
\sum_{\alpha} C(\alpha)\, \p_ju(\alpha) +
\sum_{\ui\in \hat Q_j}C'_j(\ui)[\  e_{i_2}...e_{i_{n-1}} ,e_{i_n}]_{b_j(\ui)}
= 0.\eqno(2.17)
$$ Here
$$ b_j(\ui) = a(j\ui)q_{ji_n},\eqno(2.18)
$$
 and we notice that
$$
\prod_{s\in S^c_{n-1}}b_j(s\ui) = \prod_{s \in S_n^c}a(s j\ui) = \prod
q_{i_1i_2},~j = 1,...,k,
\eqno(2.19)
$$ where one or the other (but not both unless $Q = 1^n$) of the first two
products may run more than once over an orbit of
$S^c_{n-1}$, resp. $S^c_n$. The third product is over all imbeddings of
$(i_1,i_2)$ into  $Q$. We conclude that, if the cocycle condition (2.13) holds,
then in each
$\cb_{Q_j}$, with commutation factors $b_j(\ui)$,
 all the long orbits are singular, and that if (2.13) does not hold, then none
of the orbits in
$\cb_{Q_1},...,\cb_{Q_k}$ are singular. In this latter case the first part of
(2.17) is absent.

We wish to show that  Eq.(2.16) requires that $C(\alpha) = 0$. Consider
first the case that $Q = \underline n :=\{1,...,n\}$; that is, the case that
there is no repetition in $Q$, each index occurring just once. Then all the
orbits are long and  the collection $\{u(\alpha)\} = \{u(\ui)\} =
\{e_{i_1}...e_{i_{n-1}}e_1\}$ contains exactly one monomial for each orbit of
$S^c_n$ in $\hat Q$. Hence
$$
\sum C(\alpha)u(\alpha) = ve_1,~~v\in \cb_{Q_1},
$$ and for $j \neq 1$,
$$
\sum C(\alpha) \p_ju(\alpha) = (\p_jv)e_1
$$ contains at most one monomial from each orbit of $S^c_{n-1}$. Since all
these
orbits are singular, the first term in (2.17) plays exactly the same role as
the first term in (2.15), complementing ${\bf C}\hat A_{Q_1}$ in $\cb_{Q_1}$.
In other words, the two sums in (2.17) must vanish separately, $\p_jv = 0$ for
all $j$, and hence by the stipulation, $v = 0$. The theorem is proved for this
case.

\b
\def\un{\underline{n}}

\no{\it 2.2.5. Lemma.} Let $\un  = \{1,2,...,n\}$  and let $\phi:
\un \rightarrow Q, i' \mapsto i$ be the  natural mapping that takes e.g. $1234$
to $1122$. Let $\{q'_{i'j'}  = q_{ij}\}$ be the parameters in
$\cb_{\underline n}$.  Let $C$ be a constant in $\cb_{\underline n}$, then a
constant in
$\cb_Q$ is obtained by   replacing $e_{i'}\mapsto e_i, i' = 1,...,n$.
Conversely, every constant in
$\cb_Q$ can be obtained in this way.
\b
\no{\it Proof.} The direct statement is clear. To prove the converse, let $C$
be any homogeneous constant in $\cb_Q$. In $C$, replace each generator
$e_i, i = 1,...,k$, by the sum $\sum_{i'}e_{i'}$, where the sum runs over all
$i' \in \phi^{-1}(i)$.  The resulting polynomial $C'$  is a constant.
Furthermore, each homogeneous component of $C'$
 is a constant, and one of these components, $C'_1$, say, is in
$\cb_{\underline n}$. And now $C$ is recovered  from $C_1'$ by replacing
$e_{i'} \mapsto e_i, i' = 1,...,n$.
\b The proof of the theorem is completed by noting that $C_1'$ is certainly in
${\bf C}\hat A_{\underline n}$, since there is no repetition in $\un$.
Expressing $C'_1$ as a sum of commutators (with the requisite commutation
factors determined by the parameters), we finally replace $e_{i'}$ by $e_i$ to
obtain an expression for $C$ as an element of ${\bf C}\hat A_Q$. The theorem is
proved.

 \b
\no{\it 2.2.6. Corollary.} Assume that ${\bf C} \hat A_Q$ fails to span
$\cb_Q$. Then the cocycle condition (2.13) holds and  (2.19) shows that the
long orbits (if any) in $\hat Q_1,...\hat Q_k$ are singular.  Then usually
there are constants in $\cb_Q$, the only exception being the case $Q = 1^n$,
when all the orbits of $S_{n-1}$ are short, and $q_{11} = 1$. Hence part (a) of
the stipulation is redundant if  $Q \neq 1^n$.
\b

\b
\no{\bf 3. Constants and iterated commutators.}

See definitions in 2.1.3, and let
$$
\hat X^{j} := \{X^{j   i_2...i_n }, i_2...i_n  \in \hat Q_j\},~~ j =
1,...,k~,\eqno(3.1)
$$ For parameters in general position we have $\cb_Q = {\bf C}\hat X = \cup_{j
= 1,...,k} {\bf C}\hat X^j$, a disjoint union. In general we shall show that,
not only are the constants in ${\bf C}\hat X$; they are in fact all contained
in each and every one of the subspaces ${\bf C}\hat X^j,~ j = 1,...,k$.

\b
\no{\it 3.1. Constants in ${\bf C}\hat X^j$}.

Here we shall first construct an algorithm for the determination of constants,
without at first inquiring about its effectiveness.
\b

\no{\it 3.1.1.} Define $Y_j^{i_2...i_n}$ by
$$
\partial_j X^{i_1...i_n} =
\delta^{i_1}_j(1-q_{i_1i_2}q_{i_2i_1})Y_j^{i_2...i_n},\eqno
$$  then
$$ Y_j^{i_2}= e_{i_2},~ Y_j^{i_2...i_p} =
[Y_j^{i_2...i_{p-1}},e_{i_p}]_{b_j(i_2...i_p)},
$$ with commutation factors
$$ b_j(i_2...i_p) = q_{ji_p}a(ji_2...i_p),
$$ already introduced in (2.18) for $p = n$.

Let $\hat Y_j := \{Y_j^{i_2...i_n} ;~ i_2...i_n  = \hat Q_j\}$, where $\hat
Q_j$ is the set of permutations of $Q_j$.
\b
\no{\it 3.1.2.} Fix $j \in \{1,...,k\}$. A constant in ${\bf C}\hat X^j$ is an
element
 $\sum_{{\underline i} \in \hat Q_j} C({\underline i}) X^{ji_2...i_n}$ (where
${\underline i}$ is short for
$i_2...i_n$)   such that, for $n > 2$,\footnote*{Under the stipulation, when
$n > 2$, the factor in parenthesis cannot vanish.}
$$
\sum_{{\underline i}\in \hat Q_j} C({\underline i})
(1-q_{ji_2}q_{i_2j})Y_j^{i_2...i_n} = \sum_{{\underline i}\in \hat Q_j}
C({\underline i})
(1-q_{ji_2}q_{i_2j})[Y_j^{i_2...i_{n-1}},e_{i_n}]_{b_j({\underline i})} =
0.\eqno(3.2)
$$ Consequently, a constant in ${\bf C}\hat X^j$ implies that there is a
relation among the elements of the set
$\{\hat Y_j^{i_2...i_n},~ i_2...i_n \in \hat Q_j\}$. Conversely, if such a
relation exists, then either there is a constant in ${\bf C}\hat X^j$, or else
there is a relation among the elements of the set
$\hat X^j$. By the stipulation, neither situation can arise if the set $Q$ is
replaced by $Q_i$, any
$i$, so we can conclude that, for each value of $i_n$, the set of
$Y_j^{i_2...i_{n-1}}$'s that appear inside the commutator (3.2) is linearly
independent. Therefore, the existence of a solution of (3.2) is equivalent to
the existence of a solution of the following equation,
$$\sum_{{\underline i}\in \hat Q_j } C'({\underline i}) [e_{i_2}...e_{i_{n-1}}
,e_{i_n}]_{b_j({\underline i})} = 0.\eqno(3.3)
$$ By Proposition 2.1.2, this in turn depends on the cocycle condition for the
commutation factors; namely, such relations imply that there is $s\in S_{n-1}$
such that
$$
\prod_{{\underline i} \in [sQ_j]^c}b_j({\underline i}) = 1,\eqno(3.4)
$$  where $[sQ_j]^c$ is the orbit through $sQ_j$ of the cyclic subgroup
$S_{n-1}^c$ of the group $S_{n-1}$ of permutations of $i_2...i_n$. So we have,
in all cases, the following result.
\b
\no{\it 3.1.3.  Proposition.} For every orbit $[sQ_j]^c$ of $S_{n-1}^c$ for
which the constraint (2.19) holds, there is a linear combination of the
elements in $\hat X^j$ that is annihilated by $\p_j$; and this subset of
${\bf C}\hat X^j$ spans the subspace of constants in ${\bf C}\hat X^j$.
\b
\no{\it 3.1.4. Remark. Examples.} This set of constants is
not necessarily linearly independent; it may even be empty. The simplest
example of this is $Q = 11$ and $q_{11} = 1$, which satisfies the cocycle
condition $q_{11}^2 = 1$. A less trivial example is $Q = 112$, $Q_1 = 12$, with
$q_{11}q_{12}q_{21} = -1$, a solution of the cocycle condition
$(q_{11}q_{12}q_{21})^2 = 1$; if this is the only constraint then there are no
constants in $\cb_{112}$.

In the case of an orbit of (maximal) length $n$, we found - Eq.(2.19) -~that
the condition (3.4) reduces to
$$
\prod_{{\underline i}\in [\tilde Q_j]^c}a({\underline i}) = \prod q_{ij} =
1,\eqno(3.5)
$$  where ${\underline i} = i_1...i_n,~\tilde Q_j = jsQ_j\,$     and the second
product runs over all distinct imbeddings of  $\{i,j\}$ in
$Q$.  Complications arise when
  one or the other of the orbits
 $[sQ_j]^c$ or $[ \tilde Q_j]^c$ (but not both unless $Q = 1^n$) is short. In
the first case the condition (3.5), though it is always necessary, may not be
sufficient. In the second case the constant found by solving (3.3) may turn out
to vanish identically. An example of this is the case $Q = 1122$; the cocycle
condition $q_{11} q_{22} (q_{12}q_{21})^2 = 1$ on the short orbit through
$1212$ is stronger than the cocycle condition on the long orbit through $1122$
and when it holds there are no constants. This will all be sorted out in
3.2.5.

The case when $Q = 12...n$, without repetitions, when none of the above
complications can arise, was investigated elsewhere. The result has
implications for the general case envisaged here. The stipulation 2.2.3 remains
in force.
\b
\no{\bf 3.1.5. Theorem.} [FG]  When $Q = 12...n$, without repetitions, then
there is no constant in
$\cb_Q$ unless the constraint (3.5) holds. In this case the space of constants
in $\cb_Q$  has dimension $(n-2)!$ and is a subspace of ${\bf C}\hat X^j$ for
each $j = 1,...,k$.
\b

\no{\it 3.1.6.} In this case the co-dimension of ${\bf C}\hat X$ in $\cb_Q$ is
$(n-1)! $ (the number of distinct orbits of
$S^c_n$ in $\hat Q$) and each of the $n$ sets $\hat X^j$ is linearly
independent. The $(n-1)!$ relations among the  elements of $\hat X$ relate
the constants in ${\bf C}\hat X^i$ to the constants in  ${\bf C}\hat X^j$.

In the general case we have an algorithm for constructing all the constants in
${\bf C}\hat X^j$, starting with Eq.(3.2), but it remains to be shown that
every ${\bf C}\hat X^j$  contains all the constants in $\cb_Q$. This can be
shown by using Theorem 3.1.5 and Lemma 2.2.5, but we shall present a proof that
does not rely on 3.1.5.

\bb
\no{\it 3.2. Lyndon words.}

\b
\no{\it 3.2.1. Lyndon words and good words.} We shall introduce another basis
for $\cb_Q$. The starting point is the natural basis of distinct monomials
$e_{i_1}...e_{i_n}$ that we shall now abbreviate as
$i_1...i_n$. Let us refer to the generators $e_i = i$ as `letters'.

Step One. In any monomial $i_1...i_n$ introduce parentheses as follows.
Starting at the left end, suppose the sequence descends (in the sense of the
natural order) as far as (and including) $i_a$; open a parenthesis to the left
of $i_a$. This
$i_a$ may be immediately repeated, any number of times. Once past these
repetitions, close the parenthesis just before the first occurrance of an $i_b
\leq i_a$. Then continue towards the right, applying the same rules to the rest
of the sequence. Here are some examples,
$$ (1122),~(12)(12),~(122)1, ~2(112),~ 2(12)1,~ (22)(11).
$$ Step Two. A word is either a letter or a group of letters in a parenthesis
of the type that was constructed. (That is, the sequence inside a word has the
form $(i_a)^mi_ri_s...$ where $i_r>i_a, i_s > i_a,...$.) Order the words
lexicographically. Then apply the procedure of Step One to the words. In the
example all the words but one is already in descending order and only one new
parenthesis appears,
$$ (1222),~((12)(12)),~(122)1 ,~2(12)1,~(22)(11).
$$  This process comes to an end after a finite number of steps. However, we
shall not need to go beyond the first step. We refer to the sequences obtained
after the first step as the good sequences, and to the associated words as good
words. Good words contain no nested parentheses.

Next, we interpret each word as an iterated commutator:
$$ (ij...) = (e_ie_j...) = X^{ij...}.\eqno(3.6)
$$
The set of  all
compound words that result from applying Step One to each of the permutations
of any specific monomial  form a basis for the corresponding linear space. In
particular, the set of all good compound words associated with $Q$ is a basis
for $\cb_Q$.
\b
\no{\it 3.2.2.} Now we fix $j \in \{1,...,k\}$ and notice that, in this basis,
when the parameters are in general position, the subspace $M_j$ of all compound
words annihilated by $\p_i,~ i \neq j$, consists of all those words in which
each and every parenthesis begins with $j$, and that the dimension of this
subspace is the same as that of ${\bf C}\hat X^j$.  In the special case, when
the cocycle condition is satisfied on one or more orbits, the dimension of
$M_j$ is stable, while that of ${\bf C}\hat X^j$ generally is not. Thus it is
obvious that $M_j$ will contain all the constants, but it is {\it a priori}
possible that  ${\bf C}\hat X^j$ do  not. It is important to establish that, in
fact, it does. Towards this end we have so far established the following.
\b
\no{\it 3.2.3. Lemma.} For each $j\in \{1,...,k\}$ there exists a fixed linear
subspace $M_j$ of
$\cb_Q$ with the following properties. (a)  For generic parameters; that is,
for parameters $q_{ij}$ in general position,  the set $\hat X^j$ is a basis for
$M_j$. In the exceptional case when the cocycle condition is satisfied on some
set of orbits (the ``singular" orbits), $M_j$ contains all the constants.
\b

\no{\bf 3.2.4. Theorem.} The space $ {\bf C}\hat X^j$, for any $j \in
\{1,...,k\}$ contains the space of constants in $\cb_Q$. Hence this space is
precisely ${\bf C}\hat X^i \cup {\bf C}\hat X^j$, for any pair $(i,j),~ i\neq
j$.
\b Note that the  stipulation 2.2.3 is in force.
\ve
\no{\it  Proof.}  The spaces $M_j$ are disjoint and span $\cb_Q$. Let $M'_j$ be
a complement of ${\bf C}\hat X^j$ in $M_j$. Then the union of all the $M'_j$'s
form a compliment to ${\bf C}\hat X$ in $\cb_Q$. We know, from 2.2.4, that such
a complement, no matter how it is chosen, cannot contain a constant. So $M'_j$
cannot contain a constant. Since $M_j$ contains all the constants, so does
${\bf C}\hat X^j$.
\b

\no{\it 3.2.5. The dimension of the space of constants.} Let $\chi$, resp.
$\chi_j$, be the number of distinct orbits in $Q$, resp. $Q_j$ on which the
cocycle condition holds, the singular orbits. Assume that there are no
constants in $\cb_{Q_j}, j = 1,...,k$. Then the dimension $\#(Q)$ of the space
of constants in
$\cb_Q$ is
$$
\#(Q) = \sum_{j=1}^k \chi_j - \chi.\eqno(3.7)
$$ In the case that there are no short orbits, in $Q$ or in $Q_1,...,Q_k$,
these numbers are
$$
\chi = {|\hat Q|\over n} = {(n-1)!\over n_1!...n_k!}, ~~\chi_i  =  {n_i(n-2)!
\over n_1!...n_k!},~~\#(Q) = {(n-2)!\over n_1!...n_k!}.\eqno(3.8)
$$  Note that it is only in the case that $n_j = 1$ that the singular orbits
in $\hat Q_j$ are in 1:1 correspondence with a basis for the space of
constants.
\b An example. Let $Q=1122$; there is one long orbit and one short orbit. On
the long orbits the cocycle condition is
$$
\prod q_{ij} = \bigl(q_{11}q_{22}(q_{12}q_{21})^2\bigr)^2 = 1.
$$ If this is not satisfied, then there are no constants. If it holds, then
$\chi_1 = \chi_2 = 1$. In $Q$ there is a long orbit, on which the cocycle
condition is satisfied. But there is also a short orbit, for which the
condition is
$ q_{11}q_{22}(q_{12}q_{21})^2 = 1.
$ If it holds, then $\chi = 2$ and $\#(Q) = 0$. If it does not hold; that is,
if
$ q_{11}q_{22}(q_{12}q_{21})^2 = -1,
$ then $\chi = 1$ and $\#(Q) = 1$. This example has applications to
superalgebras.
\b

\no{\it 3.2.6. A question of rigour.} Throughout this paper we have used
language that is slightly abusive. On the one hand, we refer to an algebra \cb~
(with differential structure) that is defined only once a set of parameters is
fixed. On the other hand, there are phrases like ``when $q_{12}q_{21} = 1$"
that suggests variation of the parameters. It would be more appropriate to
define the algebra as one generated by
$\{e_i\}$ and $\{q_{ij}\}$, with relations that make the latter commute with
the former and among themselves. Rather than impose constraints  on the
parameters, one would consider a family of quotients algebras obtained by
dividing by  polynomial ideals (polynomials in the $q_{ij}$ in this instance).
No concept of continuity with respect to
variation of the parameters is invoked.
\bb
\ve
\no{\bf 4. Multiple constraints.}

Consider the case $Q = 123$. If there are no constants in the spaces
$\cb_{Q_j}$ associated with the subsets $Q_j$, then there is precisely one
constant in $\cb_Q$ if and only
$\prod q_{ij} = 1$. But if $q_{12}q_{21} = 1$, then there is a constant in
$\cb_{Q_3}$, and in this case the cocycle condition for $\hat Q$ is unrelated
to the dimension of the space of constants.  This phenomenon is central to
further study of the structure of $\cb' =
\cb/{\cal I}_q$.

\bb

\no{\bf Appendix. Examples.}

For $i \neq j$, let $\sigma_{ij} := q_{ij}q_{ji}$.

\no\underbar{$Q = \{1^n\}$.} ~  All orbits are short, $a(1^n) =
(q_{11})^{n-1},\, b_1(1^{n-1}) = (q_{11})^n.$ There is a singular orbit in
$\hat Q$ if $q_{11}^{n-1} = 1$, then ${\bf C}\hat A_Q$ is empty and there is no
constant. There is a constant $(e_1)^n$ iff $(q_{11})^n = 1,\, q_{11} \neq 1$;
it is in  ${\bf C}\hat A_Q$.
\b
\no\underbar{$Q = \{12\}$.} ~ All orbits are long; $a(12)a(21) = q_{12}q_{21}
=:\sigma_{12}$. Singular case is $\sigma_{12} = 1$, then the constant is
$[e_1,e_2]_{q_{21}} \in \hat A_Q$. The constant $[e_1,e_2]_{q_{21}} = X^{12}$
is proportional to the constant $[e_2,e_1]_{q_{12}} = X^{21}$.
\b
\no\underbar{$Q = \{1^{m+1}2\}$.} ~The unique orbit in $\hat Q$ is long;
cocycle condition
$
\bigl((q_{11})^m\sigma_{12})\bigr)^{m+1} =1.
$
   The orbit in $Q_1$ is long and thus  singular. The orbit in $Q_2$ is
short,
$b_2(1...1) = (q_{11})^m\sigma_{12}$.

(a) $(q_{11})^m\sigma_{12} = 1.$ All orbits are singular; one constant, the
Serre relation
$X^{21...1}$. ($\chi = \chi_1 =\chi_2 = 1$)

(b)  $q_{11}\sigma_{12} \neq 1.$ The short orbit is not singular;  no constant.
($\chi = \chi_1 = 1, \chi_2 = 0$)

\b
\no\underbar{$Q = \{123\}$.} ~ All orbits are long, and singular when
$\sigma_{123} := \sigma_{12}\sigma_{23}\sigma_{31`} = 1.$ One constant. ($\chi
= 2, \chi_1 =\chi_2 = \chi_3 = 1$) The constant was found in [F].
\b
 \no\underbar{$Q = \{1122.\}$.} ~ There are two orbits in $\hat Q$. If
$q_{11}q_{22}\sigma^2 = 1$, they are both singular; unique orbits in
$\hat Q_1, \hat Q_2$ are also singular, no constants. ($\chi_1= \chi_2 = 1,
\,\chi = 2$) ~When
$q_{11}q_{22}\sigma^2 = -1$ the short orbit in $\hat Q$ is nonsingular;
 one constant. ($\chi = \chi_1 =\chi_2 = 1$)

\b
\no\underbar{$Q = \{1123.\}$.} ~ All orbits are long, singular; 1 constant.
($\chi = 3,\, \chi_1 = 2,\, \chi_2 = \chi_3 =1$).
\b
\no\underbar{$Q = \{1^{2m+1}22.\}$.} ~There are $m+1$ orbits in $\hat Q$, all
of them long and thus singular when
$(q_{11})^{2m(2m+1)}(q_{22})^2\sigma_{12}^{2(2m+1)} = 1$. There are $m+1$
orbits in $\hat Q_1$, one long orbit in $\hat Q_2$, so if they are all singular
there is one constant. But one of the orbits in $\hat Q_1$ is short, and if
$(q_{11})^{m(2m+1)}q_{22}\sigma_{12}^{2m+1} = -1$, then it is not singular and
then there is no constant.

\b
\no\underbar{$Q = \{1^{2m}22.\}$.} ~There are $m+1$ orbits in $\hat Q$, one of
them short. There are $m$ orbits in $\hat Q_1$ and 1 in $\hat Q_2$, all long.
So there is one constant if the short orbit in
$\hat Q$ is nonsingular; that is, if $(q_{11})^{m(2m-1)}q_{22}\sigma_{12}^{2m}
= -1$, otherwise none.
\bb
\ve

\no{\bf References.}

[B] Bokut, L.A., Kang, S.-J., Lee, K.-H. and Malcomson, P.,
``Gr\"obner-Shirshov bases for Lie superalgebras and their universal covering
algebras"; {\tt math.RT/9809024}.
\b
[F] Fr\o nsdal, C., ``Generalizations and exact deformations of quantum
groups", RIMS Publications {\bf 33} (1997) 91-149. ({\tt math.QA/9606160})
\b
[FSS] Frappat, L., Sciarrino, A. and Sorba, P., ``Dictionary on Lie
superalgebras"; {\tt hep-th/9607161}.
\b
[GF] Fr\o nsdal C. and Galindo, A., ``The ideals of free differential
algebras", J. Algebra {\bf 222} (1999) 708-746. ({\tt math.QA/9806069})
\b
[GL] Grozman, P. and Leites, D., `` Defining relations for Lie superalgebras
with Cartan matrix"; {\tt hep-th/9702073}.
\b
[Kc] Kac, V.G., ``Infinite dimensional Lie Algebras: an introduction", Progress
in Mathematics Vol. 44, Birkh\"auser, 1983.
\b
[Kh] Kharchenko, V.K. ``An existence condition for multilinear quantum
operations," J. Algebra {\bf 217} (1999) 188-228.
\b
[Ks] Kashiwara, M.,  ``On crystal basis of the Q-analogue of universal
enveloping algebras", Duke Math. J. {\bf 63} (1991) 465-516.
\b
[L] Lusztig, G., ``Introduction to quantum groups", Progress Mathematics,
Vol.110, Birkh\"auser 1993.
\b
[M] Moody, R.V. ``Lie algebras associated with general Cartan matrices,"
Bull.Am. Math.Soc. {\bf 73} (1967) 217.
\b
[R] Rosso, M., ``Quantum Groups and quantum shuffles",  Invent. Math. {\bf 133}
(1998) 399-416.
\b
[V] Varchenko,  A., ``Multidimensional hypergeometric functions and
representation theory of Lie algebras and quantum groups" Adv. Series in Math.
Phys., Vol 21, World Scientific 1995.

\ve

\end